\documentclass[a4paper,12pt]{article}
\usepackage[english]{babel}
\usepackage[T2A]{fontenc}
\usepackage[cp1251]{inputenc}
\usepackage{amsthm}
\usepackage[tbtags]{amsmath}
\usepackage{amsfonts,amssymb}
\begin{document}

\newtheorem{theorem}{Theorem}
\newtheorem{lemma}{Lemma}
\newtheorem{proposition}{Proposition}
\newtheorem{Cor}{Corollary}

\begin{center}
{\large\bf Completely Centrally Essential Rings}
\end{center}
\begin{center}
Oleg Lyubimtsev\footnote{Nizhny Novgorod State University, Nizhny Novgorod, Russia; email: oleg\_lyubimcev@mail.ru .},
Askar Tuganbaev\footnote{National Research University MPEI, Moscow, Russia; Lomonosov Moscow State University, Moscow, Russia; tuganbaev@gmail.com .}
\end{center}
\textbf{Abstract.} A ring $R$ is said to be centrally essential if for every its non-zero element $a$, there exist non-zero central elements $x$ and $y$ with $ax = y$. A ring $R$ is said to be completely centrally essential if all its factor rings are centrally essential rings. It is proved that completely centrally essential semiprimary rings are Lie nilpotent; noetherian completely centrally essential rings are strongly Lie nilpotent (in particular, every such a ring is a $PI$-ring). Every completely centrally essential ring has the classical ring of fractions which is a completely centrally essential ring. If $R$ is a commutative domain and $G$ is an arbitrary group, then any completely centrally essential group ring $RG$ is commutative.


The work of Oleg Lyubimtsev is supported by Ministry of Education and Science of the Russian Federation, project FSWR-2023-0034, and scientific and educational mathematical center "Mathematics of technologies of the future". The study of Askar Tuganbaev is supported by grant of Russian Science Foundation, project 22-11-00052, https://rscf.ru/en/project/22-11-00052.

\textbf{MSC2020 database 16D25, 16R99}

\section{Introduction}\label{sec1}

We consider only associative rings with non-zero identity element. A ring $R$ is said to be \textsf{centrally essential} if for every its non-zero element $a$, there exist non-zero central elements $x$ and $y$ with $ax = y$. \footnote{It is clear that the ring $R$ with center $Z$ is centrally essential if and only if the module $R_{Z}$ is an essential extension of the module $Z_{Z}$.} Centrally essential rings were studied in \cite{MT18}, \cite{MT19a}, \cite{MT20c}, \cite{LT24}, \cite{MT20} and other works. A ring is said to be \textsf{completely centrally essential} if all its the factor rings are centrally essential rings. 

We give an example of a non-commutative completely centrally essential ring.

\textbf{1.1. Example.} Let $F = \mathbb{Q}(x,y)$ be the field of rational functions. We consider two partial derivations $d_1 = \dfrac{\partial}{\partial x}$ and 
$d_2=\dfrac{\partial}{\partial y}$. Then the ring $R$ consisting of matrices 
$$
\left\lbrace\left.\begin{pmatrix}
f&d_1(f)&g\\
0&f&d_2(f)\\
0&0&f
\end{pmatrix}\;\;\right|\;\; f,g\in F \right\rbrace
$$
is a non-commutative centrally essential ring; see \cite[Example 3.3]{LT24}. We consider the ideal
$$
H = \left\lbrace\left.\begin{pmatrix}
0&0&g\\
0&0&0\\
0&0&0
\end{pmatrix}\;\;\right|\;\; g\in F \right\rbrace.
$$
Since $R/H \cong \mathbb{Q}(x,y)$ and $H$ is the least ideal of $R$, the ring $R$ has no non-zero proper ideals not equal to $H$. Consequently, $R$ is a non-commutative, artinian uniserial, completely centrally essential ring. 

We denote by $Z(R)$, $P(R)$ and $J(R)$ the center, the prime radical and the Jacobson radical of the ring $R$, respectively. It is well known that Jacobson radical of right artinian ring is nilpotent. The least integer $n$ with $J(R)^n = 0$ is called the \textsf{nilpotence index} of $J(R)$. A ring $R$ is said to be \textsf{local} if $R/J(R)$ is a division ring. A ring is said to be \textsf{strongly bounded} if every its non-zero one-sided ideal contains a non-zero ideal. A ring is said to be \textsf{invariant} if all its one-sided ideals are ideals. A ring $R$ is said to be \textsf{semiperfect} if the factor factor $R/J(R)$ is artinian and all idempotents of $R/J(R)$ are lifted modulo $J(R)$ to idempotents of $R$. A ring $R$ is said to be \textsf{semiprimary} if it is semiperfect and the radical $J(R)$ is nilpotent. A module $M$ is said to be \textsf{uniserial} if the set of all submodules of $M$ is linearly ordered with respect to inclusion. A ring $R$ is said to be \textsf{right uniserial} if $R$ is a uniserial right $R$-module. We recall that an element $r$ of the ring $R$ is said to be \textsf{right regular} or a \textsf{left non-zero-divisor} if $x = 0$ for every $x\in R$ with $rx = 0$. We note that one-sided zero-divisors of a centrally essential ring are two-sided; see \cite[Lemma 2.2]{MT20b}. A ring $R$ has the \textsf{right} (resp., \textsf{left}) \textsf{classical ring of fractions} $Q_{\text{cl}}(R_r)$ (resp., $Q_{\text{cl}}(R_l)$) if and only if for any two elements $a, b\in R$, where $b$ is regular, there exist elements $c, d\in R$ such that $d$ is regular and $bc = ad$ (resp., $cb = da$). These conditions are called \textsf{Ore conditions}. If the both rings $Q_{\text{cl}}(R_r)$ and $Q_{\text{cl}}(R_l)$ exist, then they are isomorphic over $R$. In this case, one says that there exists the two-sided ring of fractions $Q_{\text{cl}}(R)$.

It is well known that every associative ring $R$ can be considered as a Lie ring with respect to the Lie multiplication $[r, s] = rs - sr$ for all $r, s\in R$. For all $A, B \subseteq R$, the additive subgroup of $R$ generated by all products $[a, b]$ (where $a\in A$ and $b\in B$) is denoted by $[A, B]$. We set $R^{[1]} = R$ and for any $i\in \mathbb{N}$, $i > 1$, $R^{[i]} = [R^{[i-1]}, R]$. If there exists a positive integer $n$ with $R^{[n+1]} = 0$, then the ring $R$ is called a \textsf{Lie nilpotent} ring. The least $n$ with this property is the \textsf{Lie nilpotence class} $R$. We set $R^{(1)} = R$, $R^{(2)}$ is the ideal generated by all elements $rs - sr$ where $r, s\in R$, and, inductively, $R^{(n)}$ is the ideal generated by all $ab - ba$ where $a\in R^{(n-1)}$, $b\in R$, $n > 1$. A ring $R$ is said to be \textsf{strongly Lie nilpotent} if $R^{(m)} = 0$ for some $m$; e.g., see \cite{C00}. 

\textbf{1.2. Remark.} In \cite[Proposition 3.3]{MT18}, it is proved that semiprime centrally essential rings are commutative. Therefore, if $R$ is a completely centrally essential ring, then the factor ring $R/P(R)$ (in particular, $R/J(R)$) is commutative. 

\textbf{1.3. Remark.} The class of completely centrally essential rings is closed with respect to homomorphic images. 

The main results of the paper are Theorems 1.4, 1.5, and 1.6.

\textbf{1.4. Theorem.} Let $R$ be a completely centrally essential ring.
\begin{enumerate}
\item[\textbf{1.}] 
If $R$ is a semiprimary ring and its Jacobson radical $J(R)$ has the nilpotence index $n$, then $R$ is a Lie nilpotent ring and its nilpotence class does not exceed $n$. 
\item[\textbf{2.}] 
If $R$ is a noetherian ring, then $R$ is strongly Lie nilpotent; in particular, $R$ is a $PI$-ring.
\end{enumerate}

\textbf{1.5. Theorem.} Every completely centrally essential ring $R$ has the classical ring of fractions $Q_{\text{cl}}(R)$ which is completely centrally essential ring.

\textbf{1.6. Theorem.} Let $R$ be a commutative unital domain and let $G$ be an arbitrary group. Then the completely centrally essential group ring $RG$ is commutative.

In connection to Theorems 1.4 and 1.6, we give the following remarks.

\textbf{1.7. Remark.} There exist centrally essential rings which are not Lie nilpotent. Indeed, \cite[Theorem 1.5]{MT20c} contains an example of a centrally essential ring which is not a $PI$ ring.

\textbf{1.8. Remark.} Theorem 1.6 does not hold for arbitrary centrally essential group rings. For example, the group algebra $RG = \mathbb{Z}_2 Q_8$ is a finite non-commutative centrally essential ring, where $\mathbb{Z}_2$ is the residue field modulo $2$; see \cite{MT18}. 

\section{Proof of Theorem 1.4}

The following lemma is proved in \cite[Theorem 1.7(3)]{LT24} is important; therefore we give its proof here.

\textbf{2.1. Lemma.} Every completely centrally essential ring is invariant.

\textbf{Proof.} Let $R$ be a centrally essential ring, $I$ be a right ideal, and let $0\neq a \in I$. Since the ring $R$ is centrally essential, $0\neq ac = d\in Z(R)\cap I$. Then $dR$ is an ideal of $R$ and $dR \subseteq I$. Consequently, the ring $R$ is strongly bounded. If, in addition, $R$ is a completely centrally essential ring, then all its factor rings are also strongly bounded. In \cite[Proposition 6]{BT88}, it is proved that $R$ is an invariant ring in this case.~$\square$

\textbf{2.2. Lemma.} \cite[Proposition 3.4]{MT19a} Let $R$ be a centrally essential semiperfect ring. Then $R/J(R)$ is a finite direct product of fields; in particular, $R/J(R)$ is a commutative ring. In addition, $R$ is a finite direct product of centrally essential local rings.

For ideals $I$ and $J$ of the ring $R$ with $I\subset J$, the factor ring $J/I$ is called a \textsf{central factor} if $[J, R]\subseteq I$ or, equivalently, $J/I$ is contained in the center $Z(R/I)$ of the ring $R /I$. A chain $(J_i)$ ideals of $R$ is called a \textsf{central series} of $R$ if every  factor $J_{i+1}/J_i$ is central. In \cite{J42}, it is proved that the ring is strongly Lie nilpotent if and only if it has a finite central series.

\textbf{2.3. The completion of the proof of Theorem 1.4.} 

\textbf{1.} By Lemma 2.2, we can assume that $R$ is a local ring. We use the induction on $k$ to verify that $[r_0, r_1, \ldots, r_k]\in J^k$ for any $r_i\in R$. If $k = 1$, then $[r_0, r_1]\in J^1$, since the factor ring $R/J(R)$ is commutative by Remark 1.2. We have that 
$[r_0, r_1, \ldots, r_k] = [[r_0, r_1, \ldots, r_{k-1}], r_k]$. By the induction hypothesis, the element $j = [r_0, r_1, \ldots, r_{k-1}]$ is contained in $J^{k-1}$. If $r_k\in J$, then $[j, r_k]\in J^k$. Let $r_k$ be invertible. By Lemma 2.1, completely centrally essential rings are invariant. Therefore, $jR = Rj$ is a two-sided ideal of $R$ and $jr_k = rj$ for some $r\in R$. Since $R$ is a centrally essential ring, $0\neq jc\in Z(R)$ for some $c\in Z(R)$. Then $0\neq jr_kc = rjc$, whence $(r_k - r)jc = 0$ and $r_k - r\in J$. As a result, we obtain
$$
[j, r_k] = jr_k - r_kj = rj - r_kj = (r - r_k)j\in J^k.
$$

\textbf{2.} Since the ring $R$ is invariant, the set of nilpotent elements coincides with prime radical $P(R) = P$; see \cite[Theorem 2]{Th60}. If $P = 0$, then $R$ is a commutative ring by Remark 1.2. Let $P\neq 0$. Since the ring $R$ is noetherian, $P$ is a nilpotent ideal. If $k$ is the nilpotence index of the ideal $P$, then $P^{k-1}\cap Z(R)$ is a non-zero ideal of $R$. Indeed, since $R/P$ is commutative, the element $rs - sr$ is contained in $P$ for all $r, s\in R$. Then $x(rs - sr) = 0$ for every $x\in P^{k-1}\cap Z(R)$. Therefore, $(xr)s = x(rs) = x(sr) = s(xr)$ and $xr\in P^{k-1}\cap Z(R)$. We similarly have $rx\in P^{k-1}\cap Z(R)$. Therefore, $P^{k-1}\cap Z(R) = P_1$ is an ideal of $R$. In addition, this ideal is non-zero, since the ring $R$ is centrally essential. Then $0\subset P_1$ and $\overline{P} = P/P_1$ is a nilpotent ideal of the factor ring $\overline{R} = R/P_1$ with nilpotence index, e.g. $k_1$. Since the factor ring $\overline{R}/\overline{P}$ is commutative and the ring $\overline{R}$ is centrally essential, we have that $\overline{P}^{k_1-1}\cap Z(\overline{R}) = \overline{P}_2$ is a non-zero ideal in the ring $\overline{R}$. We denote by $P_2$ the corresponding ideal of $R$ and obtain $0\subset P_1 \subset P_2$. By induction, we obtain a chain of ideals
$$
0 = P_0 \subset P_1 \subset P_2 \subset\ldots\subset P_m = P \subset R \eqno(1)
$$
which is finite, since the ring $R$ is noetherian. In addition, $P_{i+1}/P_i\subseteq Z(R/P_i)$ for all $i = 0, 1,\ldots,m - 1$. Consequently, chain (1) is a finite central series of ideals of $R$ and the ring $R$ is strongly Lie nilpotent.

\section{Proof of Theorem 1.5}

\textbf{3.1. Proposition.} If a centrally essential ring $R$ has the right classical ring of fractions $Q_{\text{cl}}(R_r)$, then $Q_{\text{cl}}(R_r)$ is a centrally essential ring.

\textbf{Proof.} Let $0\neq ab^{-1}\in Q_{\text{cl}}(R_r)$, where $b$ is a regular element of $R$. Since $R$ is a centrally essential ring, $0\neq ac = d$ for some $c, d\in Z(R)\subseteq Z(Q_{\text{cl}}(R_r))$. Then $0\neq ab^{-1}c = b^{-1}d$ in the ring $Q_{\text{cl}}(R_r)$. Let's assume that $bd\in Z(R)$. It follows from $(bd)q = q(bd)$ that for every $q\in Q_{\text{cl}}(R_r)$, we have that $dq = b^{-1}qdb$ and $q(db^{-1}) = (db^{-1})q$. Then $b^{-1}d\in Z(Q_{\text{cl}}(R_r))$ and $Q_{\text{cl}}(R_r)$ is a centrally essential ring. If $bd\notin Z(R)$, then $0\neq (bd)c' = b(dc')\in Z(R)$. Since $dc'\in Z(Q_{\text{cl}}(R_r))$, we pass to the case considered above.~$\square$

\textbf{3.2. The completion of the proof of Theorem 1.5.} 

It follows from Lemma 2.1 that any completely centrally essential ring $R$ is invariant. Therefore, Ore conditions hold and $R$ has the classical ring of fractions $Q_{\text{cl}}(R) = Q$. It follows from Proposition 3.1 that $Q$ is a centrally essential ring. Let $0\neq I$ be an ideal of $Q$, $Q/I = \overline{Q}$ and $\overline{0}\neq \overline{ab^{-1}}\in \overline{Q}$ where $a, b\in R$. We denote by $I^c$ the restriction of the ideal $I$ to $R$. To avoid confusion, we denote by $\widetilde{r}$ residue classes of elements $r\in R$ in the ring $R/I^c$; the ring $R/I^c$ is denoted by $\widetilde{R}$. It follows from $a\notin I$ that $a\notin I^c$ and $\widetilde{a}\neq \widetilde{0}$ in the ring $\widetilde{R}$. Since $\widetilde{R}$ is a centrally essential ring, $\widetilde{0}\neq \widetilde{a}\,\widetilde{c} = \widetilde{d}$ for some $\widetilde{c}, \widetilde{d}\in Z(\widetilde{R})$. We note that if $\widetilde{x}\in Z(\widetilde{R})$, then $\overline{x}\in Z(\overline{Q})$. Indeed, for every $r\in R$, we have $[x, r]\in I^c$. Therefore, $[x, r]\in I$, i.e., $\overline{x}\,\overline{r} = \overline{r}\,\overline{x}$ in the ring $\overline{Q}$. Therefore, 
$\overline{c}, \overline{d}\in Z(\overline{Q})$ and $\overline{0}\neq \overline{ab^{-1}}\overline{c} = \overline{b^{-1}}\,\overline{d}$ in the ring $\overline{Q}$. If $\widetilde{b}\,\widetilde{d}\in Z(\widetilde{R})$, then $(\overline{b}\,\overline{d})\overline{q} = \overline{q}(\overline{b}\,\overline{d})$ for all $\overline{q}\in \overline{Q}$. Therefore, $\overline{b^{-1}}\,\overline{d}\in Z(\overline{Q})$ and $\overline{Q}$ is a centrally essential ring. If $\widetilde{b}\,\widetilde{d}\notin Z(\widetilde{R})$, then $\widetilde{0}\neq (\widetilde{b}\,\widetilde{d})\widetilde{c'}\in Z(\widetilde{R})$. In this case, $\overline{0}\neq \overline{ab^{-1}}(\overline{c}\overline{c'})\in Z(\overline{Q})$ and $\overline{Q}$ is a centrally essential ring.~$\square$

\section{Proof of Theorem 1.6}

We recall that the group $G$ is said to be \textbf{hamiltonian} if $G$ is not abelian and an arbitrary subgroup of $G$ is normal. 

\textbf{4.1. Lemma.} Let $R$ be a commutative unital domain and let $G$ be a non-abelian group. If the group ring $RG$ is completely centrally essential ring, then $G$ is a hamiltonian group.

\textbf{Proof.} By Lemma 2.1, all completely centrally essential rings are invariant. In \cite[Lemma 1(ii)]{M79}, the assertion of the lemma is proved for right invariant rings in the case, where $R$ is a field. This proof can be passed to the case of a commutative unital domain. We give the proof for the sake of completeness. 

It is sufficient to prove that all cyclic subgroups of $G$ are normal. Let $a, g\in G$ and $H = <a>$. We consider the left ideal $I = RG(1 - a)$ which is two-sided, since the ring $RG$ is invariant. Then $g^{-1}(1 - a)g\in I$ and $1 - g^{-1}ag = \alpha(1 - a)$ for some $\alpha\in RG$. We consider a mapping $\Theta: RG\to RG$ such that $\Theta(\sum_{g\in G}r_gg) = \sum_{g\in H}r_gg$. Then $\Theta$ is an $RH$-algebra homomorphism. Therefore, $1 - \Theta(g^{-1}ag) = \Theta(\alpha)(1 - a)$. Since the element $1 - a$ is not invertible, $\Theta(g^{-1}ag)\neq 0$ and $g^{-1}ag\in H$.~$\square$

\textbf{4.2. Lemma.} Let $R$ be a ring and let $G$ be a group. If the group ring $RG$ is a completely centrally essential ring, then $R$ is a completely centrally essential ring.

\textbf{Proof.} We have $R\cong RG/\omega G$, where $\omega G$ is the fundamental ideal of the ring $RG$. Since by Remark 1.3 ring $RG/\omega G$ is completely centrally essential, $R$ is a completely centrally essential ring.~$\square$

We recall that the quaternion group $Q_8$ is a group with two generators $a$, $b$ and defining relations $a^4 = 1$, $a^2 = b^2$ and $aba^{-1} = b^{-1}$. Then $Q_8 = \{e, a, a^2, b, ab, a^3, a^2b, a^3b\}$ with center $Z(Q_8) = \{e, a^2\}$.

\textbf{4.3. The completion of the proof of Theorem 1.6.} 

Let $RG$ be a completely centrally essential ring and let $G$ be a non-abelian group. Then $G$ is a hamiltonian group by Lemma 4.1. It is well known that a hamiltonian group is a direct product $G = Q_8\times A\times B$ where $A$ is an elementary abelian $2$-group and $B$ is an abelian group in which every element is of finite order. Then
$$
RG = R(Q_8\times A\times B)\cong RQ_8(A\times B),
$$
and $RQ_8$ is a completely centrally essential ring by Lemma 4.2. Since $RQ_8$ is an invariant ring, it follows from \cite[Corollary 2.3]{GL11} that $char R = 2$ or $char R =0$ (we note that the ring $R$ is a field for $char R = 2$, by \cite[Theorem 2.4(1)]{GL11}). If $char R = 0$ and $R$ is a domain, then $RQ_8$ does not have nil-ideals, see \cite[Theorem 5.1]{ZM75}; in particular, $RQ_8$ is semiprime. Since $RQ_8$ is a centrally essential ring by assumption, it is commutative by Remark 1.2. This is a contradiction. 

Let $char R = 2$. We have that
$$
Z(RQ_8) = \{\lambda_0 + \lambda_1a^2 + \lambda_2(a + a^3) + \lambda_3(b + a^2b) + \lambda_4(ab + a^3b) \, | \, \lambda_i\in R\};
$$
e.g., see \cite[Theorem 3.6.2]{MS02}. We take the ideal $I = R\widehat Q_8$ where $\widehat Q_8 = \sum_{g\in Q_8}g$. We consider the factor ring 
$\overline{R} = R/I$ which has the basis $\{\overline{e}, \overline{a}, \overline{a}^2, \overline{a}^3, \overline{b}, \overline{ab}, \overline{a^2b}\}$ and the center 
$$
Z(\overline{R}) = \{\lambda_0\overline{e} + \lambda_1\overline{a}^2 + \lambda_2(\overline{a} + \overline{a}^3) + 
\lambda_3(\overline{b} + \overline{a^2b}) + \lambda_4(\overline{e} + \overline{a} + \overline{a}^2 + \overline{a}^3 + \overline{b} + \overline{a^2b})\}.
$$
We note that if $x\in RQ_8$ and $ax\neq xa$, then $\overline{a}\,\overline{x}\neq \overline{x}\, \overline{a}$ in $\overline{R}$. Otherwise, 
$xa = ax + \alpha\widehat Q_8$ for some $0\neq\alpha\in RQ_8$. If $x = \sum\alpha_gg$, then $(\sum\alpha_gg)a = a(\sum\alpha_gg) + \alpha\widehat Q_8$. Then for some $\alpha_i\in R$, $i = 1, 2, 3, 4$, we obtain the following contradictory equality:
$$
\alpha_1b + \alpha_2ab + \alpha_3a^2b + \alpha_4a^3b = \alpha\widehat Q_8.
$$
It is similarly verified that it follows from $bx\neq xb$ that $\overline{b}\,\overline{x}\neq \overline{x}\, \overline{b}$. We take the element 
$\overline{r} = (\overline{e} + \overline{a} + \overline{b} + \overline{ab})\notin Z(\overline{R})$. When we multiply the element $\overline{r}$ by basis elements of the center $Z(\overline{R})$, we obtain either $\overline{r}$ or $\overline{0}$. Therefore, $\overline{R}$ is not a centrally essential ring. Consequently, $RQ_8$ is not a completely centrally essential ring. This is a contradiction.~$\square$

\section{Additional Examples and Remarks}

\textbf{5.1. Example.} We consider the ring
$$
R = \left\lbrace\left.\begin{pmatrix}
a&b\\
0&c
\end{pmatrix}\;\;\right|\;\; a, c\in \mathbb{Z},\, a - c\in 2\mathbb{Z},\, b\in \mathbb{Z}_4 \right\rbrace.
$$
Since 
$$
\begin{pmatrix}
1&0\\
0&3
\end{pmatrix}\begin{pmatrix}
0&1\\
0&0
\end{pmatrix} = \begin{pmatrix}
0&1\\
0&0
\end{pmatrix}\neq \begin{pmatrix}
0&3\\
0&0
\end{pmatrix} = \begin{pmatrix}
0&1\\
0&0
\end{pmatrix}\begin{pmatrix}
1&0\\
0&3
\end{pmatrix},
$$
$R$ is not commutative. We verify that $R$ is a centrally essential ring. If $a\neq 0$ or $c\neq 0$, then 
$$
\begin{pmatrix}
a&b\\
0&c
\end{pmatrix}\begin{pmatrix}
4&0\\
0&4
\end{pmatrix} = \begin{pmatrix}
4a&0\\
0&4c
\end{pmatrix}
$$
for non-zero central matrices $\begin{pmatrix}
4&0\\
0&4
\end{pmatrix}$ and $\begin{pmatrix}
4a&0\\
0&4c
\end{pmatrix}$.
If $a = c = 0$ and $b = 1$ or $b = 3$, then
$$
\begin{pmatrix}
0&b\\
0&0
\end{pmatrix}\begin{pmatrix}
2&0\\
0&2
\end{pmatrix} = \begin{pmatrix}
0&2b\\
0&0
\end{pmatrix}.
$$
It is clear that the matrix $\begin{pmatrix}
2&0\\
0&2
\end{pmatrix}$
is contained in the center of the ring $R$. For the matrix $A = \begin{pmatrix}
0&2b\\
0&0
\end{pmatrix}$, we have
$$
\begin{pmatrix}
0&2b\\
0&0
\end{pmatrix}\begin{pmatrix}
a'&b'\\
0&c'
\end{pmatrix} = \begin{pmatrix}
0&2bc'\\
0&0
\end{pmatrix} \quad \mbox{and} \quad \begin{pmatrix}
a'&b'\\
0&c'
\end{pmatrix}\begin{pmatrix}
0&2b\\
0&0
\end{pmatrix} = \begin{pmatrix}
0&2ba'\\
0&0
\end{pmatrix}.
$$
Therefore, the matrix $A$ is central if and only if $2bc' = 2ba'$ in the ring $\mathbb{Z}_4$. However, the last equality is true, since $a'$ and $c'$ are comparable modulo 2 in the ring $\mathbb{Z}$. Consequently, $R$ is a centrally essential ring. It is easy to see that the ring $R$ is completely centrally essential. By Theorem 1.5, the ring $R$ has the classical ring of fractions $Q$; in addition, $Q$ contains inverse matrices of matrices $\begin{pmatrix}
a&b\\
0&c
\end{pmatrix}$ such that $a$, $c$ are odd integers which are comparable modulo 2. In addition, the ring $Q_{\text{cl}}(R)$ also is completely centrally essential ring.

A ring $R$ is said to be \textsf{reversible} (or \textsf{commutative in zero}) if it follows from $ab = 0$ that $ba = 0$ for $a, b\in R$. Uniserial artinian rings are reversible (in particular, uniserial artinian completely centrally essential rings are reversible). Indeed, it follows from \cite[Lemma 2.1]{MT20} that every right ore left ideal of the uniserial artinian ring $R$ is a power of the Jacobson radical $J(R)$. For ideals $I$ and $K$ of $R$, it follows from $IK = 0$ that $KI = 0$, i.e., the ring $R$ is reversible.

\textbf{5.2. Example.} Let $K = \mathbb{Z}_2$ be the residue field modulo $2$ and let $T_2(K)$ be the ring of upper triangular matrices of order $2$ over the field $K$. We consider the ring
$$
R = \left\lbrace\left.\begin{pmatrix}
k&a&b\\
0&k&a\\
0&0&k
\end{pmatrix}\;\;\right|\;\; k\in K; a,b\in T_2(K) \right\rbrace.
$$
In the matrix $A\in R$, we set $k = b = 0$, $a = \begin{pmatrix}
1&0\\
0&0
\end{pmatrix}
$;
in the matrix $B\in R$, we set $k = b = 0$, $a = \begin{pmatrix}
1&1\\
0&1
\end{pmatrix}$.
Then we have $AB\neq BA$. Then $R$ is a non-commutative ring with center
$$
Z(R) = \left\lbrace\left.\begin{pmatrix}
k&k'&b\\
0&k&k'\\
0&0&k
\end{pmatrix}\;\;\right|\;\; k\in K, k'\in D_2(K), b\in T_2(K) \right\rbrace,
$$
where $D_2(K)$ is the subring of scalar matrices of $T_2(K)$. For the matrix $C\in T_2(K)$ and the identity matrix $E\in D_2(K)$, we have
$$
0\neq \begin{pmatrix}
0&C&0\\
0&0&C\\
0&0&k
\end{pmatrix} \cdot \begin{pmatrix}
0&E&0\\
0&0&E\\
0&0&0
\end{pmatrix} = \begin{pmatrix}
0&0&CE\\
0&0&0\\
0&0&0
\end{pmatrix}\in Z(R),
$$
Therefore, $R$ is a centrally essential ring. We note that the ring $R$ is subdirectly indecomposable and has the least non-zero ideal 
$
I = \begin{pmatrix}
0&0&\begin{pmatrix}
0&K\\
0&0
\end{pmatrix}\\
0&0&0\\
0&0&0
\end{pmatrix}.
$
Since the factor ring $R/I$ is commutative, the factor rings with respect to remaining ideals are commutative. Consequently, $R$ is a non-commutative completely centrally essential ring consisting of $128$ elements. In addition, the ring $R$ is not reversible. Indeed, for matrices
$A\in R:\, k = b = 0$, $a = \begin{pmatrix}
1&0\\
0&0
\end{pmatrix}
$
and $B\in R:\,k = b = 0$, $a = \begin{pmatrix}
0&1\\
0&0
\end{pmatrix}
$
we have $AB\neq 0$ but $BA = 0$. Consequently, the ring $R$ is not uniserial. Indeed, ideals 
$$
I_1 = \begin{pmatrix}
0&\begin{pmatrix}
K&0\\
0&0
\end{pmatrix}& T_2(K)\\
0&0&\begin{pmatrix}
K&0\\
0&0
\end{pmatrix}\\
0&0&0
\end{pmatrix}, \quad 
I_2 = \begin{pmatrix}
0&\begin{pmatrix}
0&0\\
0&K
\end{pmatrix}& T_2(K)\\
0&0&\begin{pmatrix}
0&0\\
0&K
\end{pmatrix}\\
0&0&0
\end{pmatrix} 
$$
are not comparable in the ring $R$.

\textbf{5.3. Example.} There exist uniserial artinian centrally essential rings which are not completely centrally essential. We consider an example constructed in \cite[Proposition 2.6]{MT20}. Let $F$ be a field with non-trivial derivation $\delta$. We consider a mapping $f\colon F\rightarrow M_4(F)$ from the field $F$ into the ring of $4\times 4$ matrices over $F$ defined by the relation
$$
\forall a\in F,\quad f(a)=\begin{pmatrix}a&0&0&0\\
0&a&0&0\\\delta(a)&0&a&0\\0&0&0&a\end{pmatrix}.
$$
It is directly verified that $f$ is a ring homomorphism.
We set 
$$
x = \begin{pmatrix}0&0&0&0\\
1&0&0&0\\0&1&0&0\\0&0&1&0\end{pmatrix}.
$$
Then the subring $R$ of the ring $M_4(F)$ generated by the set $f(F)\cup\{x\}$ is a non-commutative uniserial artinian centrally essential ring. We consider the factor ring $\overline{R} = R/Rx^2$ which is not commutative, since $\overline{x}\overline{f(a)}\neq \overline{f(a)}\overline{x}$, where 
$\delta(a)\neq 0$. In addition, the center $Z(\overline{R})$ of the ring $\overline{R}$ consists of scalar matrices with elements $a\in F$ on the main diagonal and 
$\delta(a) = 0$. It is clear that for the element $\overline{f(a)}$, where $\delta(a)\neq 0$, there does not exist a non-zero diagonal matrix $\overline{C}\in Z(\overline{R})$ such that $\overline{f(a)}\overline{C}\in Z(\overline{R})$. Consequently, the ring $\overline{R}$ is not centrally essential, whence we have that $R$ is not a completely centrally essential ring.

A ring $R$ is said to be \textsf{semicommutative} if it follows from $ab = 0$ that $aRb = {0}$ for $a, b\in R$. 

\textbf{5.4. Example.} It follows from Lemma 2.1 that completely centrally essential rings are semicommutative. This is not true for arbitrary centrally essential rings. Let $\mathbb{Z}_2$ be the residue field modulo $2$ and let $G = D_4$ be the dihedral group of order $8$ defined by generators $a,b$ and defining relations 
$a^4 = b^2 = (ab)^2 = 1$. It is easy to verify that 
$$
G = \{1,a,a^2,a^3,b,ab,a^2b,a^3b\},\;
G' = Z(G) = \langle a^2 \rangle.
$$
It follows from \cite[Proposition 2.6]{MT18} that the group algebra $\mathbb{Z}_2G$ is centrally essential. At the same time $(1 + b)^2 = 1 + b^2 = 1 + 1 = 0$ and
$$
(1 + b)a(1 + b) = a + ba + ab + bab = 1 + a^3 + ab + a^3b\neq 0.
$$

\textbf{5.5. Remark.} In \cite{MT20c}, it is shown that the ring $R[x]$ is centrally essential provided $R$ is a centrally essential ring. For every non-commutative completely centrally essential ring $R$, the ring $R[x]$ is not completely centrally essential. Indeed, if $R[x]$ is a completely centrally essential ring, then it is invariant. However, any non-commutative ring $R[x]$ is not invariant, see \cite[Lemma 3]{HHKP95}.

\textbf{5.6. Remark.} Let $R$ be a ring and let $G$ be a hamiltonian group. If the group ring $RG$ is centrally essential, then the element $2\cdot 1$ is not invertible in $R$. Indeed, it is sufficient to consider the case $RQ_8$. If $2\cdot 1$ is invertible in $R$, then 
$RQ_8\cong R(Q_8/Q'_8)\times \omega(Q_8, Q'_8)$ is a direct ring product, where $\omega(Q_8, Q'_8)$ is the ideal of $RG$ generated by the set 
$\{h - 1:\, h\in Q'_8\}$; see \cite[Proposition 3.6.7]{MS02}. The ideal $\omega(Q_8, Q'_8)$ has the basis $\{f, af, bf, abf\}$ where 
$f = 1 - e_{Q'_8} = \frac{1}{2}(1 - a^2)$, $e_{Q'_8} = \cfrac{1}{|Q'_8|}\widehat{Q'_8}$, and the center $Z = \{\alpha_0f \, \mid \, \alpha_0\in Z(R)\}$. Then, for example, for the element $af\in\omega(Q_8, Q'_8)\setminus Z$, there does not exist an element $c\in Z$ with $0\neq af\cdot c\in Z$. Therefore, 
$\omega(Q_8, Q'_8)$ is not a centrally essential ring. Consequently, $RQ_8$ is not a centrally essential ring. This is a contradiction.

\textbf{5.7. Remark.} If $R$ is a domain and $G$ is a hamiltonian group, then the group ring $RG$ is centrally essential if and only if $R$ is commutative and $char(R) = 2$. Indeed, if $R$ is a domain and $RQ_8$ centrally essential ring, then it follows from \cite[Lemma 2.1]{MT18} that the ring $R$ is centrally essential; therefore, $R$ is commutative by Remark 1.2. If $char(R)\neq 2$, then the group ring $RQ_8$ does not contain nil-ideals by \cite[Theorem 5.1 and 5.2]{ZM75}, whence $RQ_8$ is commutative, which is impossible. Conversely, if $R$ is a commutative domain of charactristic 2, then $RQ_8$ is a centrally essential ring; see \cite[Proposition 2.6]{MT18}.

\end{document}